\documentclass[preprint,12pt]{elsarticle}

\usepackage{amssymb}
\usepackage{amsmath}
\usepackage{subfigure}
\usepackage[linesnumbered]{algorithm2e}
\usepackage{comment} 

\journal{Computer Physics Communications}

\begin{document}

\begin{frontmatter}



\title{A decomposition method with minimum communication amount for parallelization of multi-dimensional FFTs}


\author[jaist,issp]{Truong Vinh Truong Duy\corref{cor1}}
\ead{duytvt@\{jaist.ac.jp,issp.u-tokyo.ac.jp\}}
\author[jaist]{Taisuke Ozaki}
\ead{t-ozaki@jaist.ac.jp}
\address[jaist]{Research Center for Simulation Science, Japan Advanced Institute of Science and Technology (JAIST), 1-1 Asahidai, Nomi, Ishikawa 923-1292, Japan}
\address[issp]{Institute for Solid State Physics, The University of Tokyo, Kashiwanoha 5-1-5, Kashiwa, Chiba 277-8581, Japan}

\cortext[cor1]{Corresponding author. Tel.: +81 761 51 1987}

\begin{abstract}
The fast Fourier transform (FFT) is undoubtedly an essential primitive that has been applied in various fields of science and engineering. In this paper, we present a decomposition method for parallelization of multi-dimensional FFTs with smallest communication amount for all ranges of the number of processes compared to previously proposed methods. This is achieved by two distinguishing features: adaptive decomposition and transpose order awareness. In the proposed method, the FFT data are decomposed based on a row-wise basis that maps the multi-dimensional data into one-dimensional data, and translates the corresponding coordinates from multi-dimensions into one-dimension so that the resultant one-dimensional data can be divided and allocated equally to the processes. As a result, differently from previous works that have the dimensions of decomposition pre-defined, our method can adaptively decompose the FFT data on the lowest possible dimensions depending on the number of processes. In addition, this row-wise decomposition provides plenty of alternatives in data transpose, and different transpose order results in different amount of communication. We identify the best transpose orders with smallest communication amounts for the 3-D, 4-D, and 5-D FFTs by analyzing all possible cases. Given both communication efficiency and scalability, our method is promising in development of highly efficient parallel packages for the FFT.

\end{abstract}


\end{frontmatter}


\section{Introduction}
\label{Introduction}


Since the pioneering work of J. W. Cooley and J. W. Tukey in 1965 \cite{citeulike:1436613}, the fast Fourier transform (FFT) has become an essential kernel for numerous fields of science and engineering, such as electronic structure calculations, digital signal processing, medical image processing, communications, astronomy, geology, optics, and so on \cite{Haynes2000130,Clarke199214,Broughton2008,Gonzales1992}. 
In those applications, the FFT is usually performed with a large data set in multiple dimensions many times, making the FFT a computationally expensive calculation. Given the importance of the FFT and widespread of massively parallel supercomputers of multi-core processors, it should come as no surprise that there have been a lot of efforts to parallelize the FFT that can be roughly categorized into two approaches: the direct approach for one-dimensional (1-D) FFTs, and the transpose approach for multi-dimensional FFTs. In the former approach, the 1-D FFT is parallelized by dividing the original problem into sub-problems recursively, which are assigned to the processes for solving in parallel \cite{6331996,Takahashi:2006:IPF:1775059.1775222,1460927}.  


In the latter approach for multi-dimensional FFTs, the FFT data are delivered to the processes so that they can have enough data to perform the sequential 1-D FFTs locally in one specific dimension in parallel, for example with the wildly popular FFTW package \cite{1386650}. The primary concern of this paper is the communication efficiency of the latter approach, as the data must be transposed, repeatedly if necessary until the data for all the dimensions have been FFT-transformed, for conducting the 1-D FFTs in other dimensions. The data transpose requires communications among the processes, and must be considered carefully, whereas the 1-D FFT is independent and communication-free. As a result, the communication efficiency heavily depends on the domain decomposition method, which should minimize the amount of communication, while not sacrificing scalability. 

The domain decomposition can be carried out in any degree from one- to $M$-dimensions for the $M$-dimensional ($M$-D) FFT, facing the trade-off between communication efficiency and scalability. The lower the degree of the domain decomposition is, the higher the communication efficiency is, and the lower the scalability is. That said, the 1-D decomposition achieves the highest communication efficiency with the lowest scalability, in contrast to the $M$-D decomposition, which carries the lowest communication efficiency with the highest scalability. 

Let us take the most popular three-dimensional (3-D) FFT as an example, where the domain decomposition method exists in three forms: one-, two-, and three-dimensions. The alphabet hereafter is used to denote the dimensions, for example, $a$ for the first dimension, $b$ for the second dimension, etc. The 1-D decomposition \cite{Haynes2000130,Dmitruk20011921}, or slab decomposition, divides the 3-D data into equal blocks of complete $ab$-planes, for instance, along the $c$-axis, and requires only one transpose step with the smallest amount of communication, but the number of processes applicable is limited to the size of one dimension. The 3-D decomposition, or volumetric decomposition, partitions the 3-D data along all three dimensions \cite{Eleftheriou2003,Eleftheriou2005,5388792,Fang2007531}, and therefore three data transpose steps are necessary to perform the 1-D FFTs along the three dimensions, making the 3-D decomposition the costliest in terms of communication amount, yet the most scalable, as the maximum number of processes is equal to the size of the 3-D data. Sandwiched between the 1-D and 3-D decompositions is the 2-D decomposition (pencil or rod decomposition) \cite{Ayala2012,takahashi2010implementation,Pekurovsky2007,Li2010} that draws smaller amount of communication than the 3-D decomposition, and offers higher scalability than the 1-D decomposition, as the number of processes applicable is now up to the size of 2-D plane.

Our work is motivated by the fact that the aforementioned domain decomposition methods usually pre-define the dimensions of decomposition, regardless of the number of processes. In particular, the 2-D decomposition partitions in two dimensions, even when the number of processes is less than the size of one dimension and it is possible to decompose in only one dimension to reduce the communication amount, because lower degree of decomposition leads to smaller amount of communication. As can be expected, the situation is worse for the 3-D decomposition, which works in three dimensions paying no attention to the number of processes. For better communication efficiency, this factor should be taken into account in the decomposition method. 

In addition, we observe that the order of data transpose has no effect upon the communication efficiency in those decomposition methods, as they treat the dimensions involved in the decomposition in the same way. For example, as mentioned above, the 2-D decomposition divides the $ab$-plane evenly, i.e., the $a$- and $b$-axes are the same and the processes are allocated approximately equal numbers of data points along these two axes. In the two subsequent transpose steps, there is no difference in doing $ac$ or $bc$ first, because they all result in the same amount of communication. However, we find that if the dimensions involved in the decomposition are handled differently, specifically the data points on the $ab$-plane are divided to the processes in ascending order of their $a$- and then $b$-coordinates, the transpose order will become a key factor in the subsequent transpose steps. Choosing the right transpose order will actually reuse a lot of data from the previous step and consequently reduce a substantial communication amount.

Another motivation is that while the decomposition method for parallel 3-D FFTs has been extensively investigated, little work has explored beyond them so far. In fact, 4-D and 5-D FFTs have various applications, and examples can be found in fields such as lattice quantum chromodynamics simulations, where a Landau gauge fixing algorithm is implemented using the 4-D FFT \cite{Cardoso2013124}, in medical image processing with a 4-D FFT-based filtering \cite{Eklund2011}, in photography \cite{Veeraraghavan:2007:DPM:1276377.1276463,citeulike:11611828}, in drug design with 5-D FFT-based protein-protein docking algorithms \cite{Ritchie01092008,Ritchie2010,Ritchie2012}, and others \cite{Kovacs:hv0008}. Hence, there is a real need to develop parallel $M$-D FFTs beyond 3-D FFTs.

In this paper, we present a decomposition method for parallelization of multi-dimensional FFTs with two distinguishing features: adaptive decomposition and transpose order awareness for achieving smallest communication amount compared to previously proposed methods. In our method, the FFT data are decomposed based on a row-wise basis that maps the $M$-D data into 1-D data, and translates the corresponding coordinates from multi-dimensions into one-dimension for equally dividing and allocating the resultant 1-D data to the processes. As a result, our method can adaptively decompose the FFT data on the lowest possible dimensions depending on the number of processes so that the communication amount can be minimized in the first place, differently from previous works that have fixed-dimensions of decomposition. In particular, the method decomposes in one dimension if the number of processes is less than or equal to the size of one dimension, in two dimensions if the number of processes is greater than the size of one dimension and less than or equal to the size of two dimensions, up to $M$-dimensions. Another unique feature of our method is the awareness of the transpose order. The row-wise decomposition provides plenty of alternatives in data transpose, and different transpose order results in different amount of communication. By analyzing all possible cases for transpose order, we find out the best transpose orders with smallest communication amounts for 3-D, 4-D, and 5-D FFTs. Finally, our method is generalized to $M$-D FFTs.

The remainder of the paper is organized as follows. Section 2 describes $M$-D FFTs. The domain decomposition method is presented in Section 3, and Section 4 gives comparison in terms of communication amount between our method and other methods. Our study is concluded in Section 5. 

\section{Multi-dimensional FFTs}

We start the discussion by introducing the $M$-D FFT by way of the 1-D FFT. The 1-D FFT transforms a 1-D data $X(j)$ of $N$ complex numbers $(X(0),X(1),...,X(N-1))$ into another 1-D data $\bar{X}(k)$ of $N$ complex numbers $(\bar{X}(0),\bar{X}(1),...,\bar{X}(N-1))$ as follows:
\begin{equation}
\bar{X}(k) = \sum_{j=0}^{N-1} \omega_N^{jk}X(j) ,  
\end{equation} 
where $\omega_N = e^{-2 \pi i/N}$, $i= \sqrt{-1}$, and the factor is dropped for simplicity.

The 1-D FFT transforms the 1-D data that has exactly one discrete variable $j$ into another 1-D data structure. Similarly, the $M$-D FFT transforms an $M$-D data $X(j_1,j_2,...,j_M)$ that has $M$ discrete variables $j_1,j_2,...,j_M$ into another $M$-D data $\bar{X}(k_1,k_2,...,k_M)$, defined by 
\begin{equation}
\bar{X}(k_1,k_2,...,k_M)= \sum_{j_1=0}^{N_1-1}  \left (  \omega_{N_1}^{j_1k_1}  \sum_{j_2=0}^{N_2-1}  \left (  \omega_{N_2}^{j_2k_2}  ... \sum_{j_M=0}^{N_M-1}  \omega_{N_M}^{j_Mk_M}X(j_1,j_2,...,j_M) \right )\right ),  
\end{equation} 
where $\omega_{N_r} = e^{-2 \pi i/N_r}$, and $1 \leq r \leq M$.

As can be seen from the above equation, the $M$-D FFT can be computed in $M$ single steps, with each performing the 1-D FFTs along one specific dimension. This is a crucial feature exploited by the domain decomposition method for parallelization of the $M$-D FFT. The first step conducts the 1-D FFTs along the last dimension, for instance, followed by a transpose operation. The second step executes the 1-D FFTs along the dimension prior to the last dimension, and so on. Lastly, the $M$th step carries out the 1-D FFTs along the first dimension, ending the $M$-D FFT. 

\section{Domain decomposition method}
\label{Methods}
In this section, we present our domain decomposition method, starting with the 3-D FFT for ease of understanding. We then provide a general description of the method for the $M$-D FFT, and finally we describe the method for the 4-D and 5-D FFTs, and beyond them. 

\subsection{3-D FFTs}
Here we illustrate and examine our decomposition method for the 3-D FFT. Assume that the numbers of data points along the $a$-, $b$-, and $c$-axes are $N_{\mathrm{1}}$, $N_{\mathrm{2}}$, and $N_{\mathrm{3}}$, respectively, the number of processes is $N_p$, and $myid$ is the process identification. 

\begin{figure}[htb]
\centering
\includegraphics[scale=0.7,trim=0cm 0cm 0cm 0cm]{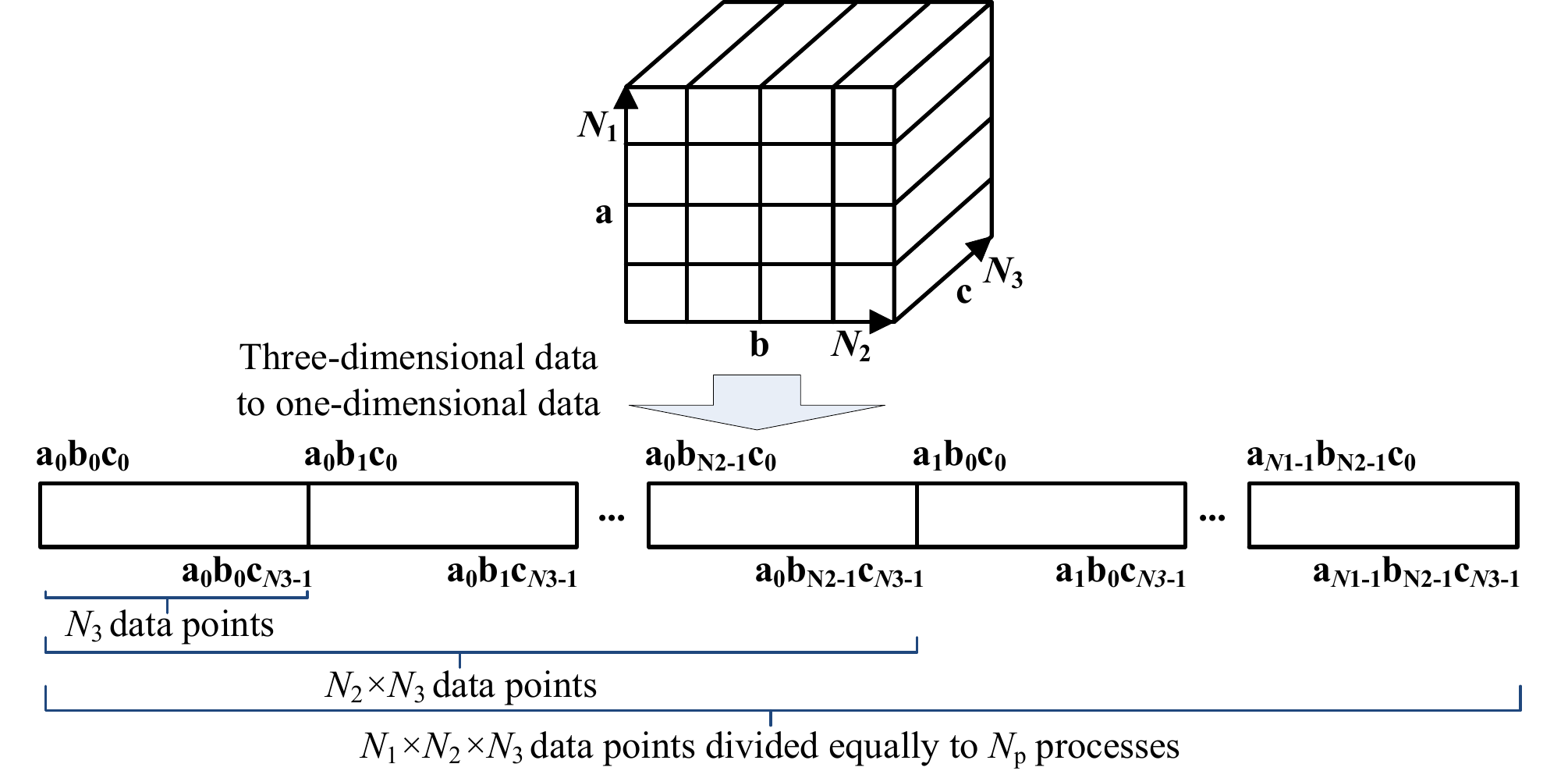}
\caption{3-D FFTs: 3-D data to 1-D data mapping in row-wise decomposition for the $abc$ order.}
\label{fig-3D-1D}
\end{figure}

Our method is basically based on a row-wise decomposition. It first translates the original 3-D FFT data $X_{\mathrm{3D}}(a,b,c)$ into a 1-D data $X_{\mathrm{1D}}(x)$, as illustrated in Fig. \ref{fig-3D-1D} for the $abc$ order as an example in a total of $3 \times 2 \times 1 = 6$ orders. 
The relationship between a 3-D coordinate $X_{\mathrm{3D}}(a_1,b_1,c_1)$ and a 1-D coordinate $X_{\mathrm{1D}}(x_1)$ in the $abc$ order is given by
\begin{equation}
x_1=a_1 \times N_2 \times N_3 + b_1 \times N_3 + c_1.
\end{equation}
It is worth mentioning that the relationships in other orders different from $abc$ can be derived in the same way. 

In order to identify the starting and ending points allocated to each process in the domain decomposition determined by $N_p$, $myid$, $N_{\mathrm{1}}$, $N_{\mathrm{2}}$, and $N_{\mathrm{3}}$, we define a function $f_{\mathrm{3D}}()$ for translating the 3-D and 1-D coordinates as follows:   
\begin{equation}
f_{\mathrm{3D}}(N_1,N_2,N_3,N_p,myid) = \begin{cases} \left \lfloor \frac{N_1 \times myid}{N_p} \right \rfloor \times N_2N_3, & \mbox{if } N_p \leq N_1 ;
\\ 
\\ \left \lfloor \frac{N_1N_2 \times myid}{N_p} \right \rfloor \times N_3 , & \mbox{if } N_1 < N_p \leq N_1N_2 ;
\\
\\ \left \lfloor \frac{N_1N_2N_3 \times myid}{N_p} \right \rfloor  , & \mbox{if } N_1N_2 < N_p \leq N_1N_2N_3,
\end{cases}
\end{equation}
where $\left \lfloor  \right \rfloor$ is the floor function. 
Our method then equally divides the resultant 1-D data to the processes, in which a process with $myid$ is assigned the data points from $X_{\mathrm{1D}}(x^{s}_{myid})$ to $X_{\mathrm{1D}}(x^{e}_{myid})$ in one dimension, where
\begin{equation}
x^{s}_{myid}=f_{\mathrm{3D}}(N_1,N_2,N_3,N_p,myid),
\end{equation}
\begin{equation}
x^{e}_{myid}=f_{\mathrm{3D}}(N_1,N_2,N_3,N_p,myid+1)-1.
\end{equation}
These 1-D coordinates can be translated back to the 3-D ones to obtain the corresponding starting and ending coordinates in three dimensions as  
\begin{equation}
\label{eq-3D-ase}
a^{(s,e)}_{myid}=\left \lfloor \frac{x^{(s,e)}_{myid}}{N_2N_3} \right \rfloor ,
\end{equation}
\begin{equation}
\label{eq-3D-bse}
b^{(s,e)}_{myid}=\left \lfloor \frac{x^{(s,e)}_{myid}-a^{(s,e)}_{myid}N_2N_3}{N_3} \right \rfloor ,
\end{equation}
\begin{equation}
\label{eq-3D-cse}
c^{(s,e)}_{myid}=x^{(s,e)}_{myid}-a^{(s,e)}_{myid}N_2N_3-b^{(s,e)}_{myid}N_3,
\end{equation}
where $X_{\mathrm{3D}}(a^{s}_{myid},b^{s}_{myid},c^{s}_{myid})$ and $X_{\mathrm{3D}}(a^{e}_{myid},b^{e}_{myid},c^{e}_{myid})$ are the starting and ending points, respectively, in three dimensions for a process with $myid$.

Consequently, the decomposition has three forms depending on the number of processes. The distribution of data points is carried out in 1-D, 2-D, or 3-D data defined by the first one, two, or three dimensions, respectively. For example, with $N_p$ processes and $N_1 < N_p \leq N_1N_2$, the decomposition in the $abc$ order takes place in the 2-D decomposition, where the data points on the $ab$-plane with the $c$-axis are divided over the processes in ascending order of their $a$- and $b$-coordinates so that each process has approximately the same number of data points on the $ab$-plane. The data points extending along the $c$-axis that have the same $a$- and $b$-coordinates are also assigned to the same process. Therefore with the 2-D decomposition on the $ab$-plane in the $abc$ order, a process with $myid$ will be assigned the data points from $X_{\mathrm{3D}}(a^{s}_{myid},b^{s}_{myid}, 0)$ to $X_{\mathrm{3D}}(a^{e}_{myid},b^{e}_{myid}, N_{\mathrm{3}}-1)$ in ascending order of the $a$-, $b$-, and $c$-coordinates, where $a^{s}_{myid}$, $b^{s}_{myid}$, $a^{e}_{myid}$, and $b^{e}_{myid}$ can be obtained from Eqs. \ref{eq-3D-ase} and \ref{eq-3D-bse}, with
\begin{equation}
x^{s}_{myid}=\left \lfloor \frac{N_1N_2 \times myid}{N_p} \right \rfloor \times N_3,
\end{equation}
\begin{equation}
x^{e}_{myid}=\left \lfloor \frac{N_1N_2 \times (myid+1)}{N_p} \right \rfloor \times N_3-1.
\end{equation}
In subsequent transpose steps, the decomposition can occur in any order of combination of the three dimensions $a$, $b$, and $c$. 

\begin{figure}[htbp]
\begin{center}
\subfigure[Transpose-order awareness: transpose from $cab$ to $cba$, a majority of data can be reused.]{\label{fig-3D-2Da}
\includegraphics[scale=0.75]{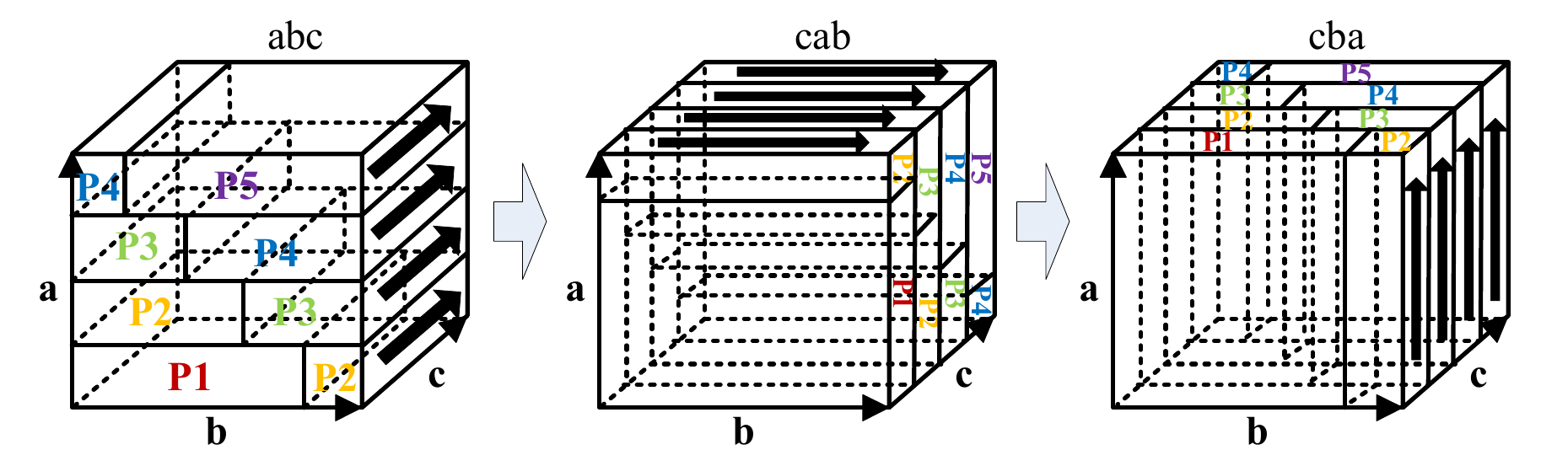}}
\subfigure[Transpose-order unawareness: transpose from $cab$ to $bca$, only a minority of data can be reused. The amount of communication is doubled compared to (a).]{\label{fig-3D-2Db}
\includegraphics[scale=0.75]{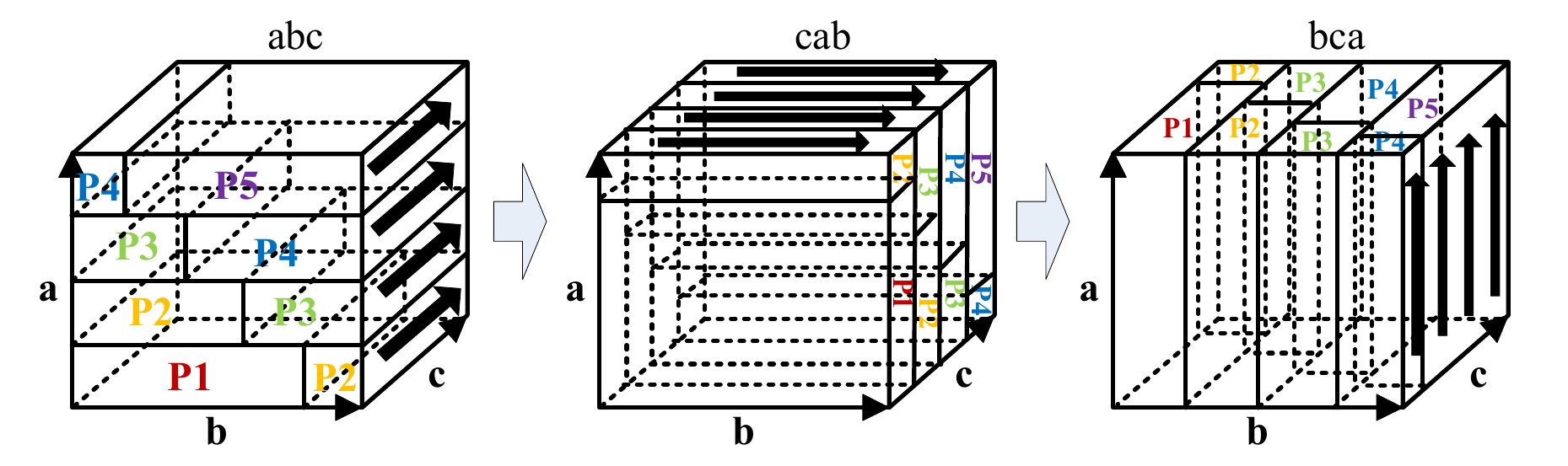}} 
\subfigure[Conventional 2-D decomposition: all dimensions are the same.]{\label{fig-3D-2Dc}
\includegraphics[scale=0.75]{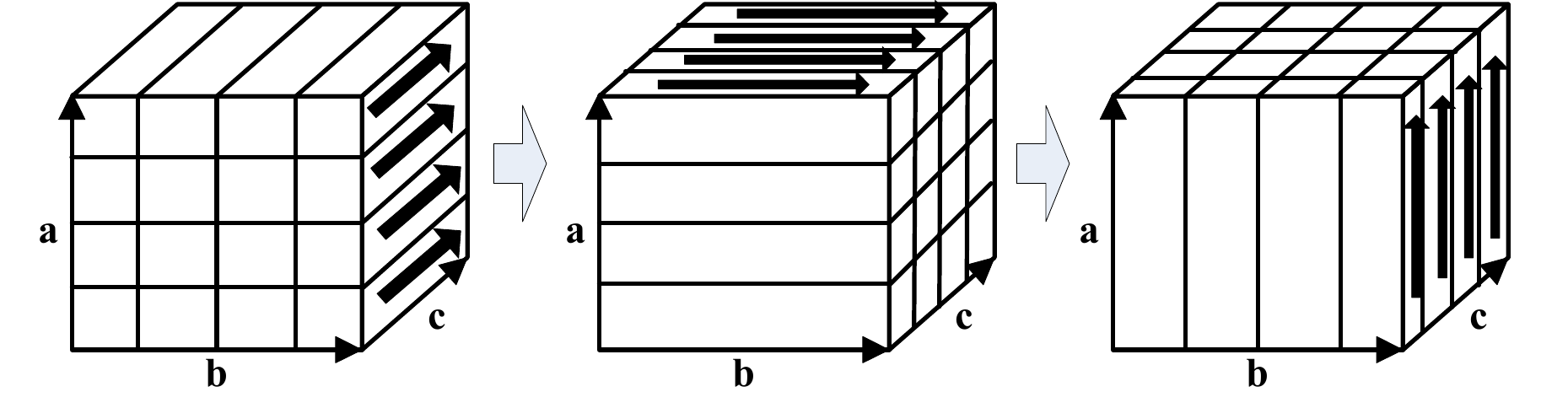}} 
\end{center}
\caption{3-D FFTs with 2-D decomposition: our method with transpose-order awareness (a) and unawareness (b), and regular 2-D method (c).}
\label{fig-3D-2D}
\end{figure}

\begin{figure}[htbp]
\begin{center}
\subfigure[8 transpose order cases and the amount of communication corresponding to the number of processes.]{\label{fig-3D-tree}
\includegraphics[scale=0.7]{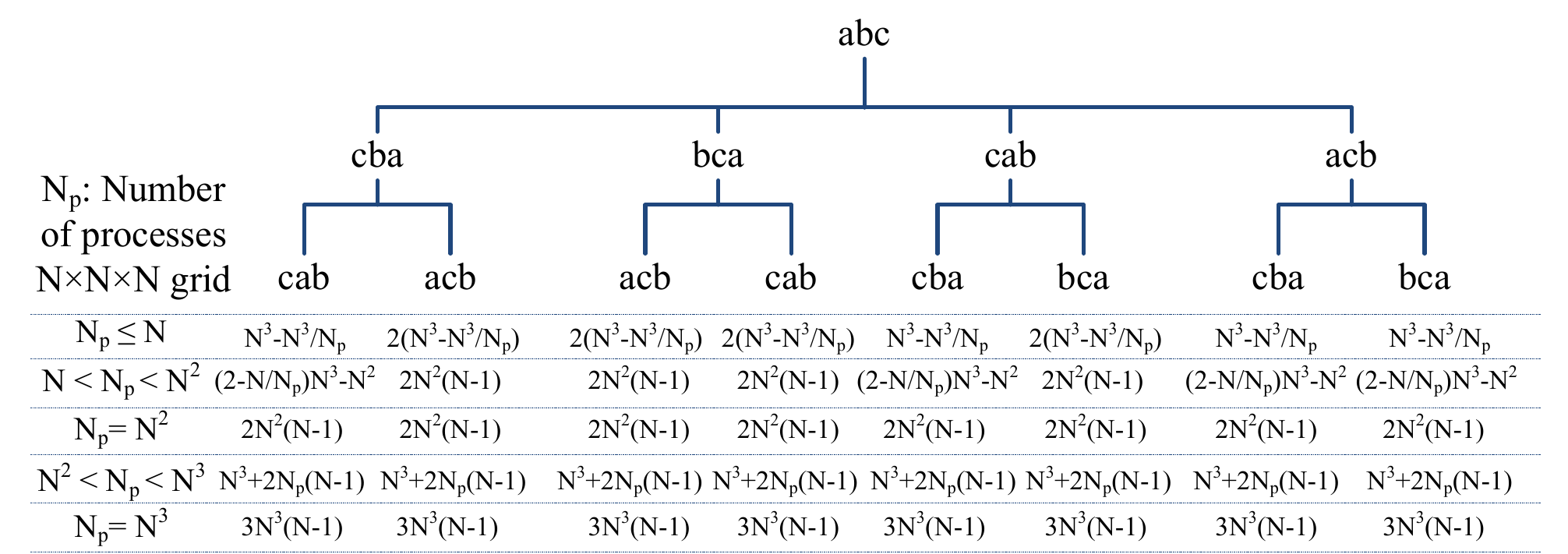}}
\subfigure[2 patterns (P1 and P2), the amount of communication corresponding to the number of processes, and the difference between their amount.]{\label{fig-3D-sum}
\includegraphics[scale=0.8]{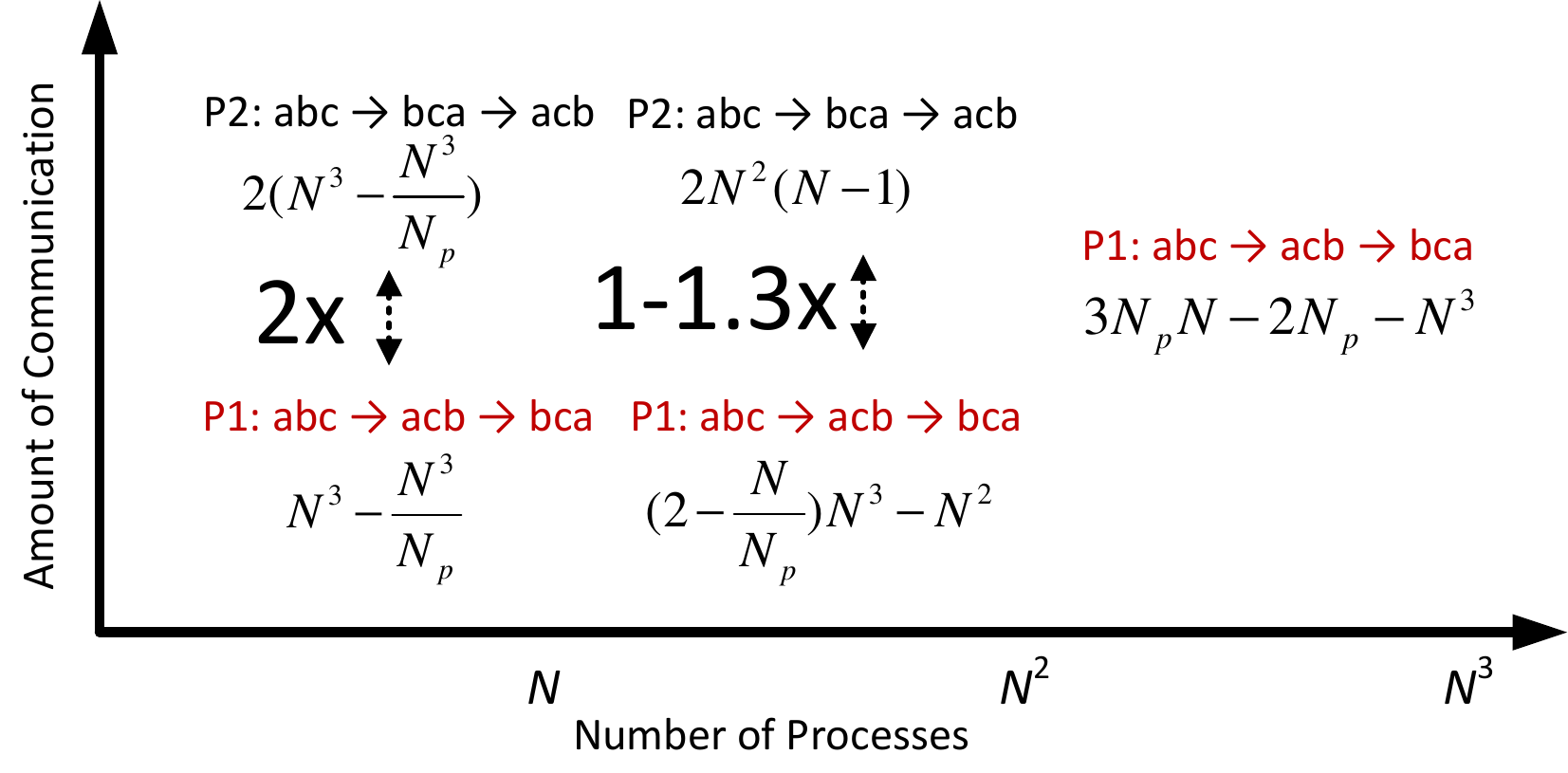}}
\end{center}
\caption{3-D FFTs: all 8 cases (a), categorized into 2 patterns (b).}
\label{fig-3D-FFT}
\end{figure}

Figure \ref{fig-3D-2D} exemplifies the operation of our method with transpose-order awareness (a) and transpose-order unawareness (b), followed by a conventional 2-D method for comparison (c). The transpose order in Fig. \ref{fig-3D-2Da} is $abc \rightarrow cab \rightarrow cba$, and in Fig. \ref{fig-3D-2Db} $abc \rightarrow cab \rightarrow bca$. Even though the only difference between them is the last decomposition, $cba$ as against $bca$, this has far-reaching implications for the amount of reused data, because a majority of data can be reused with the transpose from $cab$ to $cba$, while the transpose from $cab$ to $bca$ leaves only a minority of data that can be reused. For instance, with process P1, Fig. \ref{fig-3D-2Da} shows a large overlap between the areas assigned to it in $cab$ and $cba$, implying that a large amount of data can be reused, whereas the overlap between $cab$ and $bca$ is small (Fig. \ref{fig-3D-2Db}). In fact, as revealed later in Fig. \ref{fig-3D-tree}, the amount of communication with transpose-order unawareness, $abc \rightarrow cab \rightarrow bca$, is doubled compared to transpose-order awareness, $abc \rightarrow cab \rightarrow cba$. On the other hand, the conventional 2-D decomposition is applied to two dimensions that are treated in the same way so that the processes have approximately equal numbers of data points on these two dimensions, leaving no difference between $cba$ and $bca$, and eventually no effect of the transpose order. Also, though the illustrations are intended for 2-D decomposition, the extension to 1-D and 3-D decompositions is straightforward. 

Figure \ref{fig-3D-tree} shows a full analysis on the amount of communication to the number of processes with all 8 possible cases of transpose order for the 3-D FFT, with $N \times N \times N$ data points for the sake of simplicity and clear presentation. Interestingly, the amount of communication is grouped into only two patterns of communication in our method, and there is one pattern actually up to twice better than the other, as shown in Fig. \ref{fig-3D-sum}. The difference between these two patterns is twice, when the number of processes is up to the size of one dimension, and up to 1.3 times, when the number of processes is in between the sizes of one and two dimensions. The difference remains almost unchanged by $N$. We find four transpose orders leading to the pattern with the smaller amount of communication for all ranges of the number of processes:
\begin{displaymath}
abc \rightarrow  acb \rightarrow  bca,
\end{displaymath}
\begin{displaymath}
abc \rightarrow  acb \rightarrow  cba,
\end{displaymath}
\begin{displaymath}
abc \rightarrow  cba \rightarrow  cab,
\end{displaymath}
\begin{displaymath}
abc \rightarrow  cab \rightarrow  cba. 
\end{displaymath}
In general, we can follow any of these four orders in parallelization of the 3-D FFT. In fact, our previous work has employed the method with the 2-D decomposition and the transpose order of $abc \rightarrow  cab \rightarrow  cba$ \cite{Duy2012}. 

\subsection{General description of the method}
Let $X_{\textrm{M-D}}(N_1,N_2,...,N_M)$ be the $M$-D input data. Each process is allocated approximately $ \frac{N_1  \times N_2  \times ...  \times N_M}{N_p} $ data points. The decomposition starts from the first dimension and gradually moves down to the last dimension of the data structure to assign the data points to the processes according to the number of processes. The problem turns to finding the starting and ending coordinates of each dimension of the data structure for each process, specifically, the starting point $X_{\textrm{M-D}}(x^{s}_{1,myid},x^{s}_{2,myid},...,x^{s}_{M,myid})$ and the ending point $X_{\textrm{M-D}}(x^{e}_{1,myid},x^{e}_{2,myid},...,x^{e}_{M,myid})$ that together define the data points distributed to a particular process with $myid$. This procedure can be generalized from the previous case of 3-D FFTs.

After the $M$-D input data have been translated into the 1-D data, a 1-D coordinate $X_{\textrm{1D}}(X_1)$ is identified from the equivalent $M$-D coordinate $X_{\textrm{M-D}}(x_1,x_2,..,x_M)$ as follows:
\begin{equation}
X_1=x_1 \times N_2 \times N_3 \times \cdots \times N_M + x_2 \times N_3 \times N_4 \times \cdots  \times N_M  + \cdots  + x_{M-1} \times N_M + x_M.
\end{equation}
As well as the case of the 3-D FFT, a function $f_{\textrm{M-D}}()$ for translating the $M$-D and 1-D coordinates is defined as   
\begin{align}
f_{\textrm{M-D}}(N_1,N_2,...,N_M,N_p,myid)  = \phantom{-------------------} \nonumber \\ 
\begin{cases} \left \lfloor \frac{N_1 \times myid}{N_p} \right \rfloor \times N_2N_3 \times \cdots \times N_M, & \mbox{if } N_p \leq N_1  ;
\\
\\ \left \lfloor \frac{N_1N_2 \times myid}{N_p} \right \rfloor \times N_3 \times \cdots \times N_M , & \mbox{if } N_1 < N_p \leq N_1N_2 ;
\\
\\\cdots 
\\
\\ \left \lfloor \frac{N_1N_2 \times \cdots \times N_M \times myid}{N_p} \right \rfloor  , & \mbox{if } N_1N_2 \times \cdots \times N_{M-1} < N_p \leq N_1N_2 \times \cdots \times N_{M}. 
\end{cases}
\end{align}   
In one-dimension, the data points from $X_{\textrm{1D}}(x^{s}_{myid})$ to $X_{\textrm{1D}}(x^{e}_{myid})$ are allocated to a process with $myid$, where
\begin{equation}
x^{s}_{myid}=f_{\textrm{M-D}}(N_1,N_2,...,N_M,N_p,myid),
\end{equation}
\begin{equation}
x^{e}_{myid}=f_{\textrm{M-D}}(N_1,N_2,...,N_M,N_p,myid+1)-1.
\end{equation}
Finally, in $M$-dimensions, the starting point $X_{\textrm{M-D}}(x^{s}_{1,myid},x^{s}_{2,myid},...,x^{s}_{M,myid})$ and the ending point $X_{\textrm{M-D}}(x^{e}_{1,myid},x^{e}_{2,myid},...,x^{e}_{M,myid})$ of a process with $myid$ are given by  
\begin{equation}
x^{(s,e)}_{1,myid}=\left \lfloor \frac{x^{(s,e)}_{myid}}{N_2N_3\times \cdots \times N_M} \right \rfloor ,
\end{equation}
\begin{equation}
x^{(s,e)}_{2,myid}=\left \lfloor \frac{x^{(s,e)}_{myid}-x^{(s,e)}_{1,myid}N_2N_3\times \cdots \times N_M}{N_3\times \cdots \times N_M} \right \rfloor ,
\end{equation}
\begin{center}
...
\end{center}
\begin{align}
x^{(s,e)}_{M-1,myid}= \phantom{-------------------} \nonumber \\
\left \lfloor \frac{x^{(s,e)}_{myid}-x^{(s,e)}_{1,myid}N_2N_3\times \cdots \times N_M - x^{(s,e)}_{2,myid}N_3\times \cdots \times N_M - \cdots - x^{(s,e)}_{M-2,myid}N_{M-1}N_M }{N_M} \right \rfloor ,
\end{align}
\begin{equation}
x^{(s,e)}_{M,myid}=x^{(s,e)}_{myid}-x^{(s,e)}_{1,myid}N_2N_3\times \cdots \times N_M-x^{(s,e)}_{2,myid}N_3\times \cdots \times N_M - \cdots -  x^{(s,e)}_{M-1,myid}\times N_M. 
\end{equation}

With the definition of the function $f_{\textrm{M-D}}()$, the row-wise decomposition realizes the first feature of our method, adaptive decomposition. This is how our method is adaptive and flexible to the number of processes $N_p$. When $N_p$ is less than or equal to the size of the first dimension $N_1$, the method will partition the $M$-D data in only that dimension, and assign about $\frac{N_1}{N_p} \times (N_2 \times ... \times N_M)$ data points to each process. Likewise, if $N_p$ is between $N_1$ and $N_1 \times N_2$, the decomposition will take place in the first and second dimensions, and distribute approximately $\frac{N_1 \times N_2}{N_p} \times (N_3 \times ... \times N_M)$ data points to each process. As lower degree of decomposition requires smaller amount of communication, our method is able to minimize the communication amount in the first place. 

This adaptive decomposition is conducted when data transpose is performed to re-allocate the data points to the processes. The computation is relatively simple, and therefore its computational time is trivial. Figure 4 outlines the parallel $M$-D FFT with the decomposition for data transpose. The order of transpose used in the figure is general, and the actual order, accompanied by the amount of communication, is analyzed for the 3-D, 4-D, and 5-D FFTs. 

\begin{figure}[htb]
\label{fig.par}
\begin{algorithm}[H]
 \SetAlgoLined
 \KwIn{$X(N_1,N_2,...,N_M)$, $N_p$, $myid$}
 \KwOut{$\bar{X}(N_1,N_2,...,N_M)$}
 Step 1: Perform $N_1 \times N_2 \times \cdots \times N_{M-1}$ 1-D FFTs, each with $N_M$ points along the $M$th dimension.
 $\bar{X}^1(j_1,j_2,...,j_{M-1},k_M) = \sum_{j_M=0}^{N_M-1}  \omega_{N_M}^{j_Mk_M}X(j_1,j_2,...,j_M)$
 
 Transpose: Conduct the adaptive decomposition. 
 $\bar{X}^2(k_M,j_1,j_2,...,j_{M-1}) = \bar{X}^1(j_1,j_2,...j_{M-1},k_M)$
 
 Step 2: Perform $N_M \times N_1 \times N_2 \times \cdots \times N_{M-2}$ 1-D FFTs, each with $N_{M-1}$ points along the $(M-1)$th dimension.
 $\bar{X}^3(k_M,j_1,j_2,...,j_{M-2},k_{M-1})= \sum_{j_{M-1}=0}^{N_{M-1}-1}  \omega_{N_{M-1}}^{j_{M-1}k_{M-1}}\bar{X}^2(k_M,j_1,j_2,...,j_{M-1})$ 

 Transpose: Conduct the adaptive decomposition. 
 $\bar{X}^4(k_M,k_{M-1},j_1,j_2,...,j_{M-2}) = \bar{X}^3(k_M,j_1,j_2,...,j_{M-2},k_{M-1})$  
 
 ...
 
 Step $M$: Perform $N_M \times N_{M-1} \times \cdots \times N_2$ 1-D FFTs, each with $N_{1}$ points along the first dimension.  $\bar{X}^{2M-1}(k_M,k_{M-1},...,k_2,k_1) = \sum_{j_{1}=0}^{N_{1}-1}  \omega_{N_{1}}^{j_{1}k_{1}}\bar{X}^{2M-2}(k_M,k_{M-1},...,k_2,j_1)$ 
 
 Transpose: Conduct the adaptive decomposition.  
 $\bar{X}(k_1,k_2,...,k_{M-1},k_M) = \bar{X}^{2M-1}(k_M,k_{M-1},...,k_2,k_1)$
 
\end{algorithm}
\caption{Parallel $M$-D FFTs with adaptive decomposition for data transpose.}
\end{figure}

Since this is a row-wise-based distribution, the method leaves a number of order options to be explored, because in this case, $abc$ is no longer identical to $cab$ with the 3-D FFT, as illustrated previously. The order of transpose plays a key role in the ability of re-using data that is directly related to the amount of communication. Let us calculate the total number of cases for the $M$-D FFT. Since we start with the first transpose, it has only a single case. In the next transpose, as we have $M-1$ remaining dimensions to choose from and each dimension has $(M-1)!$ cases to be explored, there are $(M-1) \times (M-1)!$ cases in this step. Similarly, there are $(M-2) \times (M-1)!$ cases in the third transpose, because we have $M-2$ remaining dimensions with each dimension having $(M-1)!$ cases. In the $M$th (final) transpose, there are $1 \times (M-1)!$ cases. Multiplication of all the cases in $M$ steps gives us the total number of cases for the $M$-D FFT
\begin{align}
C_{\textrm{M-D}} &= \underbrace{1}_\text{1st transpose} \times \underbrace{(M-1) \times (M-1)!}_\text{2nd transpose} \times \underbrace{(M-2) \times (M-1)!}_\text{3rd transpose}  \times ... \nonumber \\ 
&\times \underbrace{2 \times (M-1)!}_\text{(\textit{M}-1)th transpose} \times \underbrace{1 \times (M-1)!}_\text{\textit{M}th transpose}  \nonumber \\ 
&= (M-1)\times(M-2)\times...\times 2 \times 1 \times (M-1)!^{M-1} \nonumber \\ 
&= (M-1)!^M 
\end{align}   

Table 1 shows the number of order cases corresponding to the number of dimensions. With the 3-D FFT, there are 8 cases only, which can be examined manually. However, with the 4-D and 5-D FFTs, these numbers are 1,296 and 7,962,624, respectively, that can be investigated thoroughly by computer simulations. Even so, it is practically difficult with 2,985,984,000,000 cases for the 6-D FFT. Therefore, we have no choice but to limit ourselves to the 5-D FFT, and try to generalize the results to those beyond them. 

\begin{table}
\label{tab-num-case}
\begin{center}
\caption{Number of cases.}
\begin{tabular}{c r}

Number of dimensions              & Number of cases \\
\hline
2 & 1  \\
3 & 8  \\
4 & 1,296  \\
5 & 7,962,624  \\
6 & 2,985,984,000,000 \\
$M$ & $(M-1)!^M$ \\
\end{tabular}
\end{center}
\label{tab-numcase}
\end{table}

\subsection{4-D FFTs, 5-D FFTs, and beyond}
\begin{figure}[htbp]
\begin{center}
\subfigure[4-D FFTs: the patterns, the amount of communication corresponding to the number of processes, and the difference between their amount.]{\label{fig-4D-sum}
\includegraphics[scale=0.7]{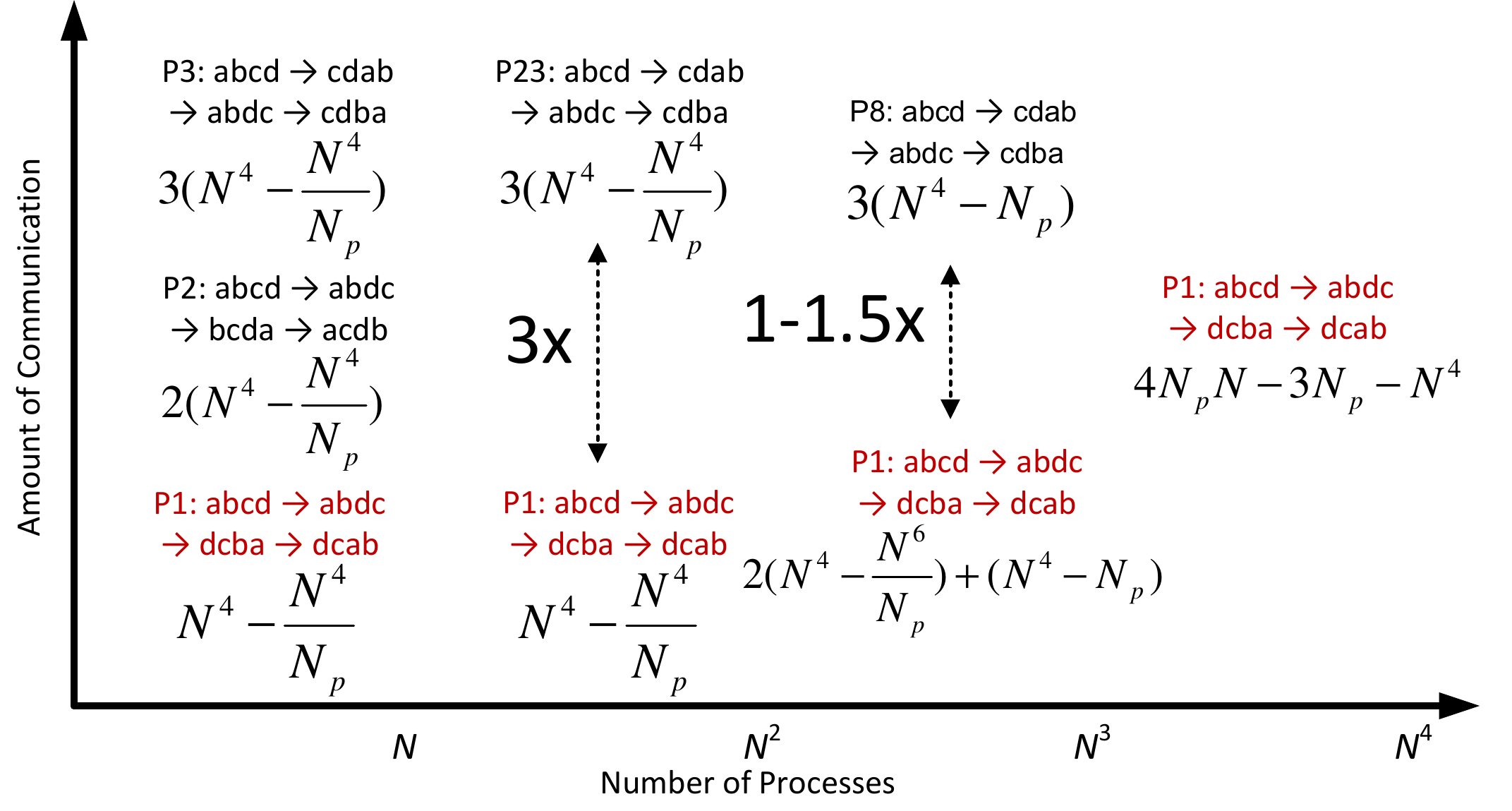}}
\subfigure[5-D FFTs: the patterns, the amount of communication corresponding to the number of processes, and the difference between their amount.]{\label{fig-5D-sum}
\includegraphics[scale=0.55]{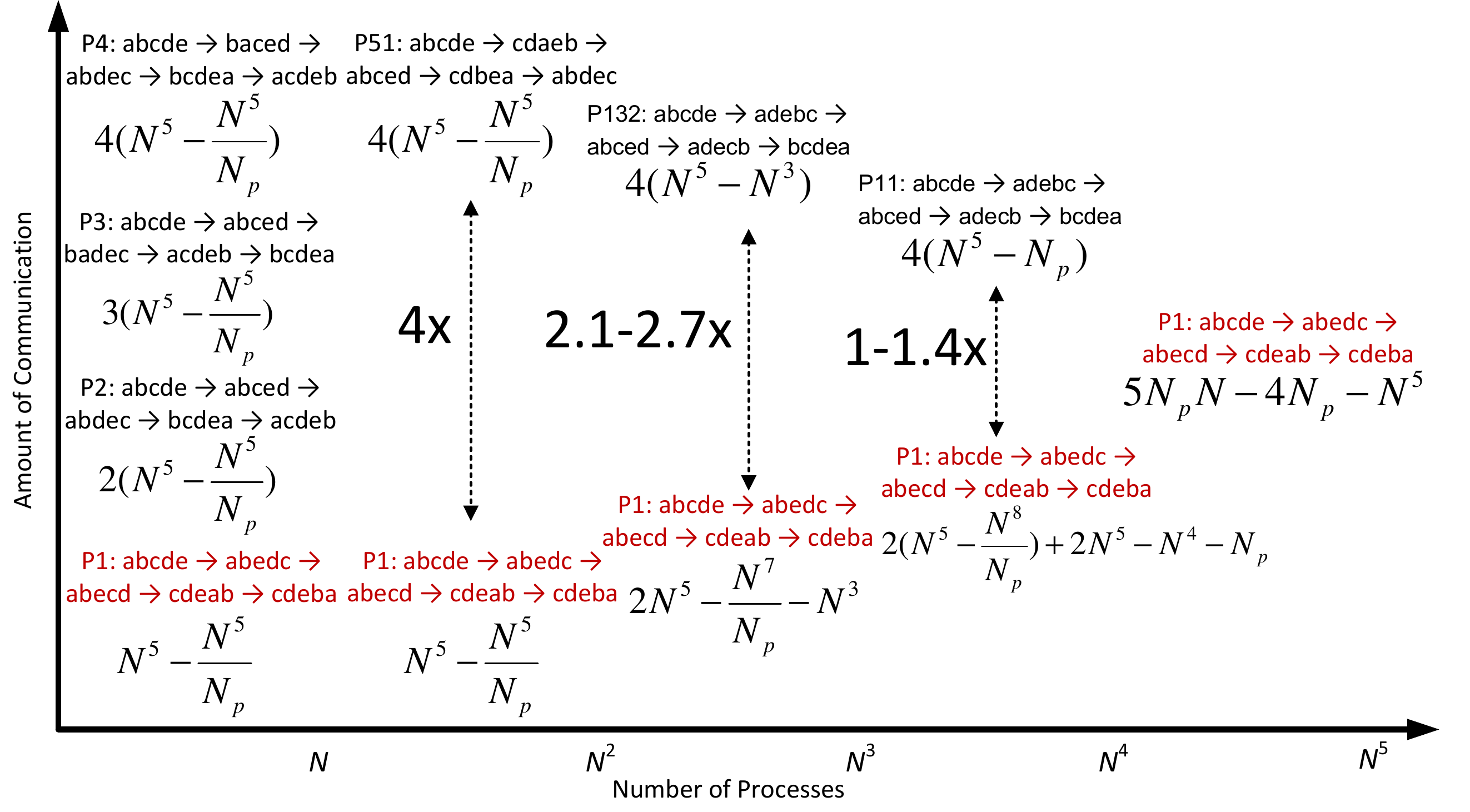}}
\end{center}
\caption{4-D and 5-D FFTs.}
\label{fig-4D-5D}
\end{figure}

Figure \ref{fig-4D-5D} depicts the patterns, the amount of communication corresponding to the number of processes, and the difference between their amount in different ranges for the 4-D FFT (a) and 5-D FFT (b), where the $M$-D data are assumed to have an equal size in all dimensions. Due to their complex nature, we do not have any illustrations for the decomposition method. As it is also difficult to present all possible cases, 1,296 for the 4-D FFT and 7,962,624 for the 5-D FFT, here we only show the best and worst patterns.  

With the 4-D FFT, we find four transpose orders: 
\begin{displaymath}
abcd \rightarrow abdc \rightarrow dcba \rightarrow dcab,
\end{displaymath}
\begin{displaymath}
abcd \rightarrow abdc \rightarrow dcab \rightarrow dcba,
\end{displaymath}
\begin{displaymath}
abcd \rightarrow abdc \rightarrow cdab \rightarrow cdba,
\end{displaymath}
\begin{displaymath}
abcd \rightarrow abdc \rightarrow cdba \rightarrow cdab,
\end{displaymath}
that always offer the smallest amount of communication for all ranges of the number of processes. The best transpose order is exactly 3 times better than the worst order for up to $N^2$ processes, and up to 1.5 times for the next range of $[N^2,N^3]$. The difference remains almost the same, apparently unaffected by $N$.  

The number of always-best transpose orders is found to be higher with the 5-D FFT than with the 4-D FFT: 96 orders, including 
\begin{displaymath}
abcde \rightarrow abced \rightarrow abedc \rightarrow cedba \rightarrow cedab,
\end{displaymath}
\begin{displaymath}
abcde \rightarrow abced \rightarrow abedc \rightarrow ecdba \rightarrow ecdab,
\end{displaymath}
\begin{displaymath}
abcde \rightarrow abced \rightarrow abedc \rightarrow cdeba \rightarrow cdeab,
\end{displaymath}
\begin{displaymath}
abcde \rightarrow abced \rightarrow abedc \rightarrow cdeab \rightarrow cdeba.
\end{displaymath}
The gap between the best and worst orders is also higher with the 5-D FFT, especially with a small number of processes. With up to $N^2$ processes, it is exactly 4 times. The gap is from 2.1 to 2.7 times for the range of $[N^2,N^3]$, and up to 1.4 times for the next range of $[N^3,N^4]$. Similar to the 3-D and 4-D FFTs, the gap is found to be almost completely independent of $N$. 

Based on these observations, our decomposition method for higher dimensional FFTs is thought to produce higher number of transpose orders that are always best for every range of the number of processes. And so is the gap in the amount of communication between the worst and best orders. Regarding the transpose order, we notice that the best orders of the 3-D FFT are in the form of either 1+2 or 2+1, i.e., $a(bc)$ or $(ab)c$, where 2 simply means a transpose $ab \rightarrow ba$. The form of $a(bc)$ leads to the subsequent transpose order of $a(cb)$, and then $(cb)a$, being one of the four best orders for the 3-D FFT. Likewise, the form of $(ab)c$ and its consequent orders of $c(ab)$ and $c(ba)$ are another best order. For the 4-D FFT, they follow the form of 2+2, $(ab)(cd)$, but neither 1+3, $a(bcd)$, nor 3+1, $(abc)d$. Meanwhile, the forms of 2+3, $(ab)(cde)$, and 3+2, $(abc)(de)$, are the best combinations for the 5-D FFT. Consequently, we expect the form of 3+3, $(abc)(def)$, to deliver better, if not best, performance for the 6-D FFT. Better forms for higher-dimensional FFTs can be derived in the same way such that their two parts are the most balanced.


\section{Comparison of Communication Amount}
\label{Comparison}
\begin{figure}[htbp]
\begin{center}
\subfigure[Our method and the 1-D method \cite{Dmitruk20011921}, 1.5-D method \cite{takahashi2012}, and 2-D method \cite{takahashi2010implementation}.]{\label{fig-3D-compare-a}
\includegraphics[scale=0.95]{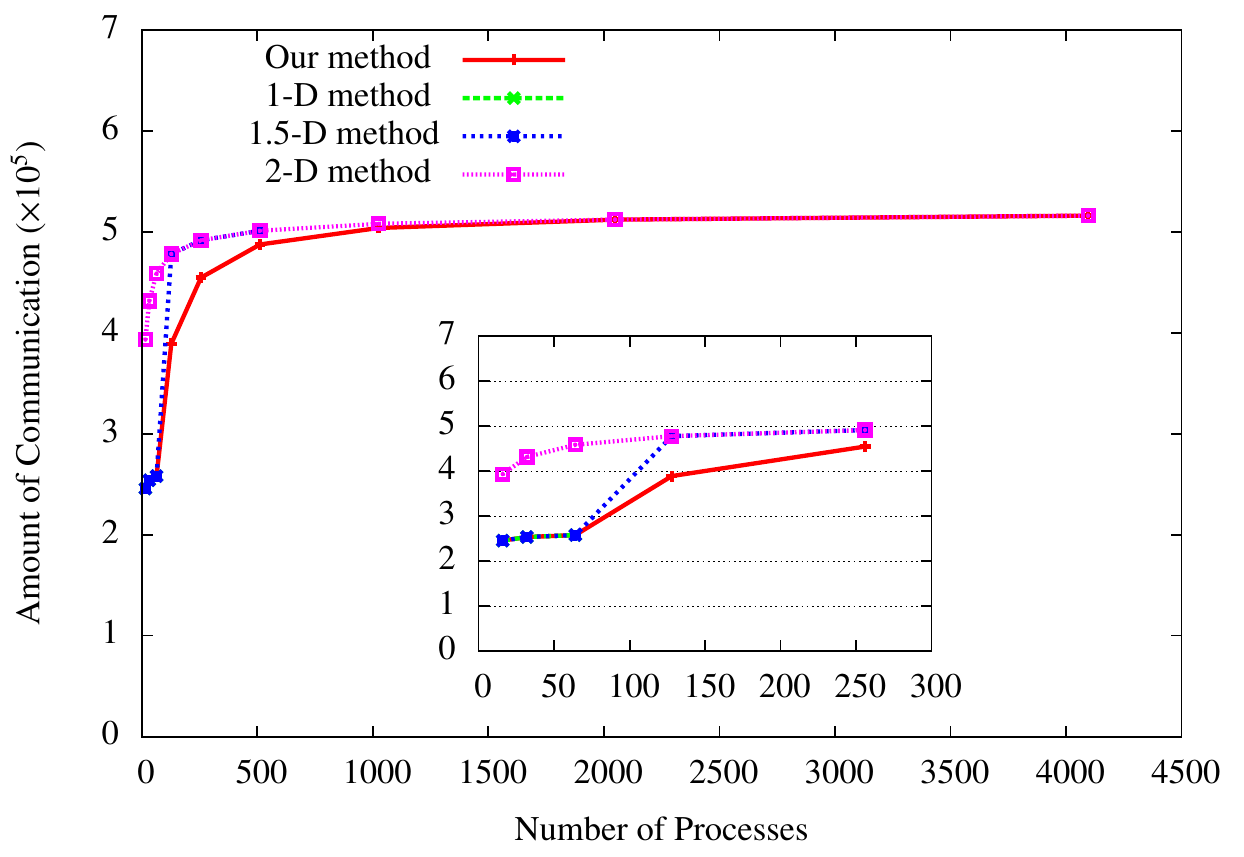}}
\subfigure[Our method and the 3-D method \cite{Eleftheriou2003}.]{\label{fig-3D-compare-b}
\includegraphics[scale=0.95]{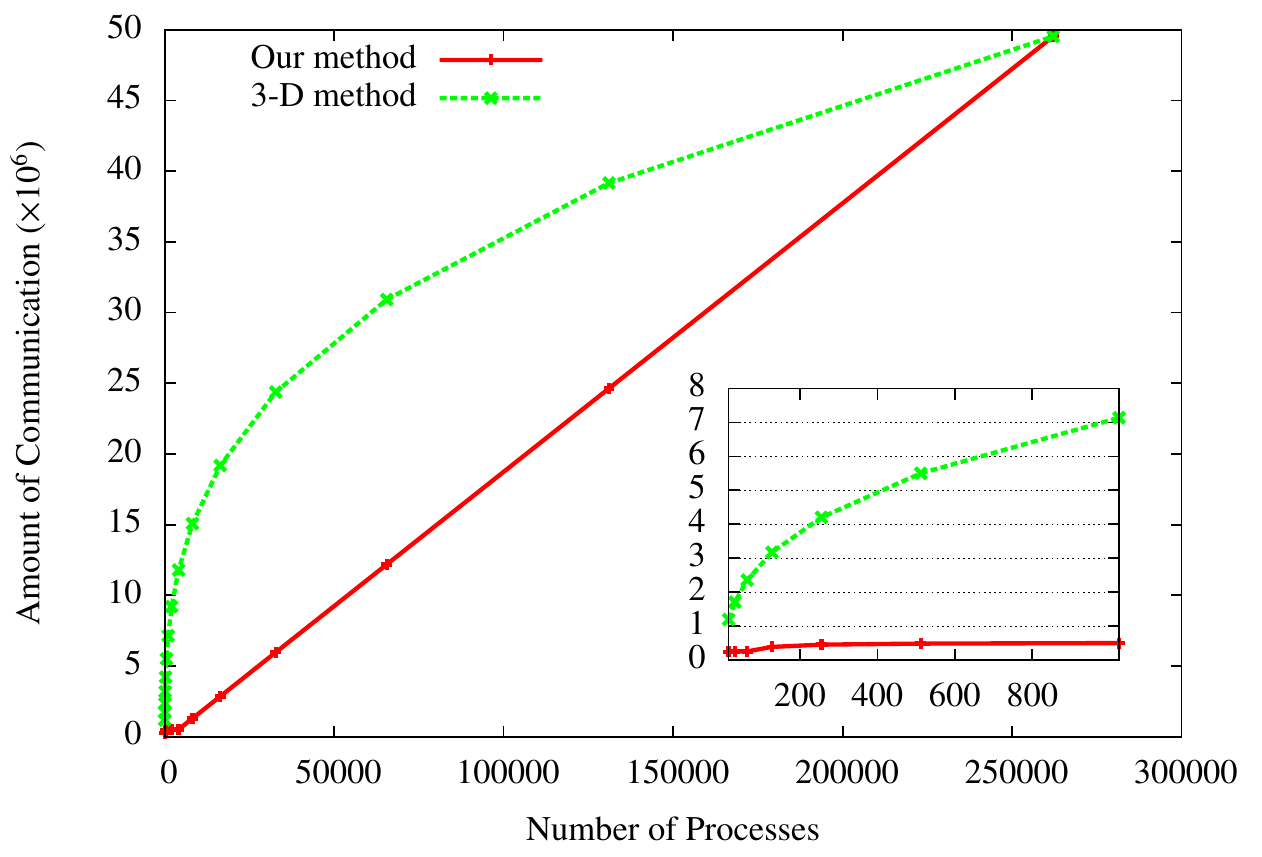}}
\end{center}
\caption{Comparison for 3-D FFTs in terms of communication amount with $64 \times 64 \times 64$ data points.}
\label{fig-3D-compare}
\end{figure}

Figure \ref{fig-3D-compare} compares our proposed method with the 1-D method \cite{Dmitruk20011921}, 1.5-D method \cite{takahashi2012}, and 2-D method \cite{takahashi2010implementation} (a), and with the 3-D method \cite{Eleftheriou2003} (b), for the 3-D FFT. The 3-D method is displayed in a separate figure for a clear presentation, as the amount of communication in this case is far larger than the other cases. The amount of communication is calculated in case of $64 \times 64 \times 64$ data points. Certainly, similar results can be produced for smaller and larger numbers of data points. 

The 1-D method divides the 3-D data into blocks of equal numbers of complete $ab$-planes along the $c$-axis, and then allocates them to the processes to perform the 1-D FFTs along the $a$- and $b$-axes, followed by a data transpose so that each process contains complete $ac$-planes to carry out the 1-D FFTs along the remaining $c$-axis, with a total communication amount
\begin{equation}
A_{1\mathrm{D}}=N^3-\frac{N^3}{N_p}.
\end{equation}
In the 2-D method, the $ab$-plane is evenly divided to the processes, with each having all the data along the remaining $c$-axis. The 1-D FFTs along $c$ are executed, followed by two transpose steps so that the 1-D FFTs along the $a$- and $b$-axes can also be performed. The 2-D method has a communication amount
\begin{equation}
A_{2\mathrm{D}}= 2N^3-2N_p(\frac{N^2}{N_p})^{\frac{3}{2}}. 
\end{equation}
The 1.5-D method lies between these two methods, with the amount being similar to that of the 1-D method, when $N_p \leq N$, and the 2-D method, when $N < N_p \leq N^{1.5}$. Lastly, the 3-D method partitions the 3-D data along all three dimensions, and requires an amount 
\begin{equation}
A_{3\mathrm{D}}= 3N^3(\frac{N}{\sqrt[3]{\frac{N^3}{N_p}}}-1). 
\end{equation}

As shown in Fig. \ref{fig-3D-compare-a}, the 1-D method is able to work only along one dimension and is limited to 64 processes, while the 2-D method decomposes the domain in two dimensions, even for fewer than 64 processes. The 1.5-D method offers a compromise between the 1-D and 2-D methods. By contrast, our method is the most flexible and adaptive, as it partitions only along one dimension when the number of processes is up to 64 on condition that it is a divisor of 64, and decomposes in two dimensions while still starting from one dimension for a larger number of processes with transpose-order awareness to reuse as many data points as possible. As a result, up to 64 processes, our method works in the same fashion as the 1-D and 1.5-D methods provided that 64 is a multiple of the number of processes, and is about 60.0\% to 77.8\% better than the 2-D method with $64 \times 64 \times 64$ data points. From this point to 512 processes, the limit of the 1.5-D method, our method still has the edge over the 1.5-D and 2-D methods. Beyond the point of 512 processes, the 1.5-D method is no longer applicable, while the two other methods can operate until reaching the limit of $64 \times 64 = 4,096$ processes for the 2-D decomposition. The performance gap also gradually decreases, however, and eventually becomes 0 from 2,048 processes. From 4,096 processes to the maximum 262,144 processes, only the 3-D decomposition is workable. Figure \ref{fig-3D-compare-b} demonstrates that our method outperforms the 3-D method, with the performance gain able to reach several orders of up to 11.6 (at 8,192 processes) for the 3-D decomposition alone. 

\begin{figure}[htbp]
\begin{center}
\subfigure[4-D FFTs with $16 \times 16 \times 16 \times 16$ data points.]{\label{fig-4D-compare}
\includegraphics[scale=0.95]{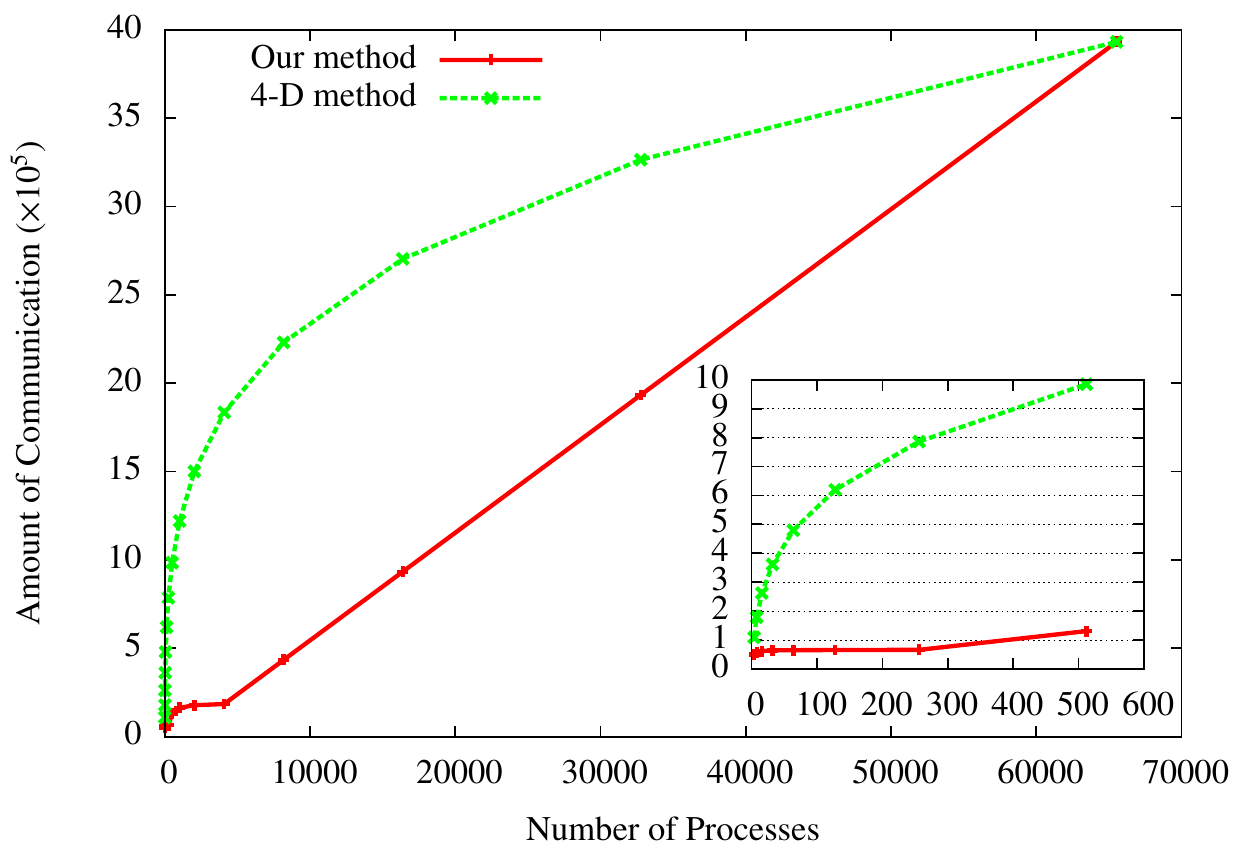}}
\subfigure[5-D FFTs with $16 \times 16 \times 16 \times 16 \times 16$ data points.]{\label{fig-5D-compare}
\includegraphics[scale=0.95]{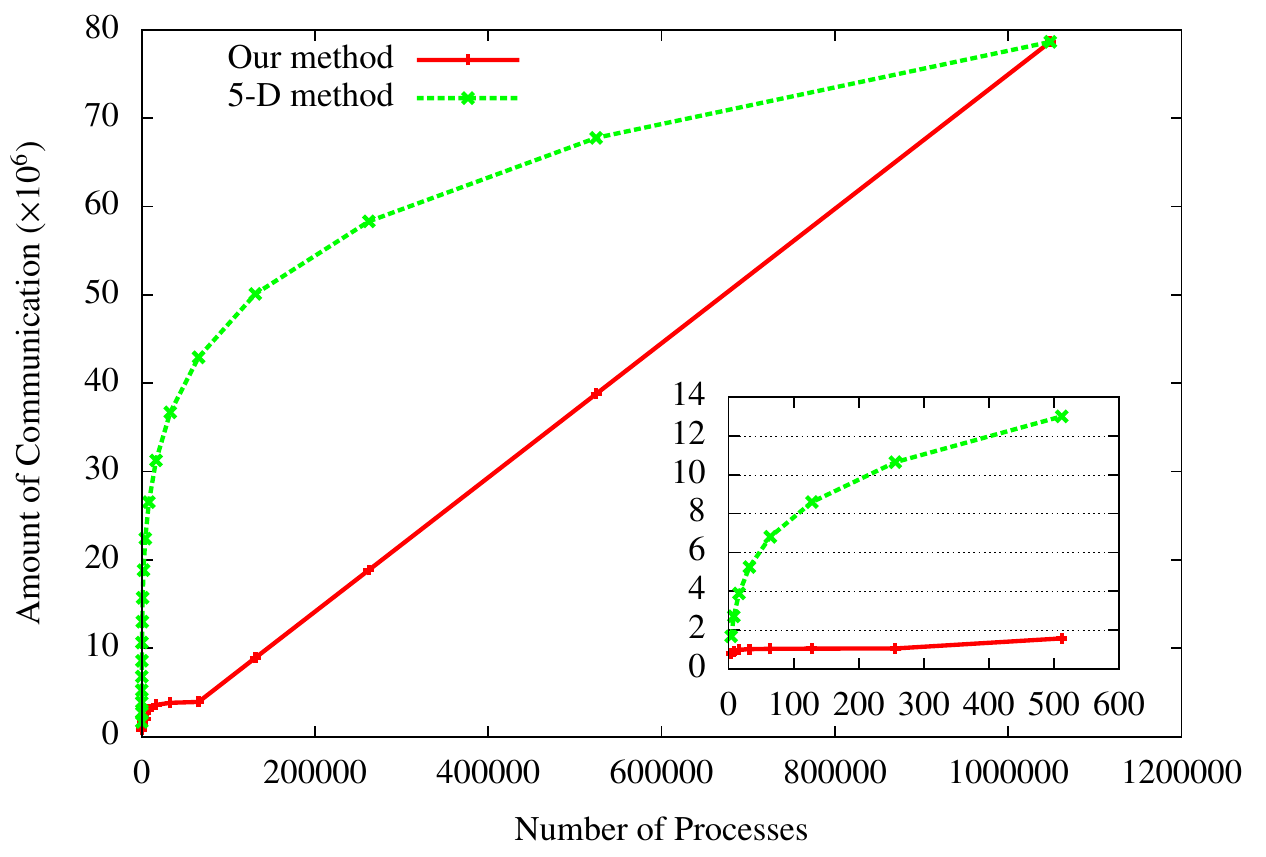}}
\end{center}
\caption{Comparison for 4-D and 5-D FFTs.}
\label{fig-4.5D-compare}
\end{figure}

Figure \ref{fig-4.5D-compare} extends the comparison by projecting the results obtained by our method and two other methods for the 4-D and 5-D FFTs, with $16 \times 16 \times 16 \times 16$ and $16 \times 16 \times 16 \times 16 \times 16$ data points, respectively, for a smaller number of processes. The 4-D and 5-D methods are assumed to operate in 4-D and 5-D decompositions, respectively, akin to the 3-D method. As a result, the amounts of communication are 
\begin{equation}
A_{4\mathrm{D}}= 4N^4(\frac{N}{\sqrt[4]{\frac{N^4}{N_p}}}-1)  
\end{equation}
for the 4-D method, and 
\begin{equation}
A_{5\mathrm{D}}= 5N^5(\frac{N}{\sqrt[5]{\frac{N^5}{N_p}}}-1)   
\end{equation}
for the 5-D method. Given these conditions, our method has a distinct advantage, leaving a wide performance gap of up to approximately 12 times for the 4-D FFT, and 11.1 times for the 5-D FFT. The gap is found to remain almost unchanged by $N$.

\section{Conclusion}
\label{Conclusion}
We have presented our decomposition method for parallelization of multi-dimensional FFTs. The communication amount is the smallest compared to previously proposed methods, accomplished by adaptive decomposition and transpose order awareness. Featured by the row-wise decomposition that translates the $M$-D data into 1-D data and evenly divides the resultant 1-D data to the processes, our method can adaptively decompose the FFT data on the lowest possible dimensions based on the number of processes. In addition, our row-wise decomposition method provides a lot of alternatives in data transpose, among them the best communication efficient orders are identified and applied. We have determined the best transpose orders for the 3-D, 4-D, and 5-D FFTs, dependent on which we find out the way for deriving the transpose orders that can deliver better performance for higher-dimensional FFTs. Comparison in terms of communication amount shows that our method is superior to other methods for the 3-D FFT, and it is anticipated to have a distinct advantage for higher-dimensional FFTs. Actually, the method has been employed in our open-source density functional theory code called OpenMX \cite{openmx}.
Boosting communication efficiency while not sacrificing scalability, our method is promising to be harnessed in development of highly efficient parallel packages for multi-dimensional FFTs.




\section*{Acknowledgements}
This work was supported by the Strategic Programs for Innovative Research (SPIRE), MEXT, and the Computational Materials Science Initiative (CMSI), and Materials Design through Computics: Complex Correlation and Non-Equilibrium Dynamics A Grant in Aid for Scientific Research on Innovative Areas, MEXT, Japan. 



\bibliographystyle{model1a-num-names}
\bibliography{biblio}







\end{document}